\documentclass[10pt]{article}
\usepackage{latexsym,amsfonts,amssymb,amsmath,amsthm,graphicx}
        \newtheorem{theorem}{Theorem}
        \newtheorem{definition}{Definition}
        \newtheorem{proposition}{Proposition}[section]
        \newtheorem{corollary}{Corollary}[section]
        \newtheorem{lemma}{Lemma}

\begin{document}

\title{The Algebra P_n is Koszul}
\author{David Nacin}

\abstract{The algebras $Q_n$ describe the relationship between the
roots and coefficients of a non-commutative polynomial. I.Gelfand,
S.Gelfand, and V. Retakh have defined quotients of these algebras
corresponding to graphs. In this work we find the Hilbert series
of the class of algebras corresponding to the graph $K_3$. We also
show this algebra is Koszul.}

\section{Koszul Algebras}

There are a number of equivalent definitions of Koszul algebras
including this lattice definition from Ufnarovskij \cite{U}.

\begin{definition} \label{AVR}

A quadratic algebra $A = \{V,R\}$ (where $V$ is the span of the
generators and $R$ the span of the generating relations in
$V\otimes V$) is Koszul if the collection of $n-1$ subspaces
$\{V^{\otimes i-1} \otimes R \otimes V^{\otimes n-i-1}\}_i$
generates a distributive lattice in $V^{\otimes n}$ for any $n$.

\end{definition}

In \cite{Musti} the following criterion is given for
distributivity of a modular lattice:

\begin{theorem} \label{MUST}

Suppose  $\{x_1, \cdots, x_n\}$ generates the modular lattice
$\Omega$. If any proper subset of $\{x_1, \cdots, x_n\}$ generates
a distributive sublattice then $\Omega$ is distributive iff for
any $2 \leq k \leq n-1$ the triple $x_1 \vee \cdots \vee x_{k-1},
x_k, x_{k+1} \wedge \cdots \wedge x_n$ is distributive.

\end{theorem}

Applying theorem \ref{MUST} to definition \ref{AVR} we get the
following corollary we will use in various chapters throughout
this work:

\begin{corollary} \label{mainkoszul}

The quadratic algebra $A = \{V,R\}$ (where $V$ is the span of the
generators and $R$ the span of the generating relations in
$V\otimes V$) is Koszul if $RV^{n - 2} \cap VRV^{n - 2} \cap
\cdots \cap V^{a - 2}RV^{n - a}$, $V^{a - 1}RV^{n - a - 1}$,
$V^{a}RV^{n - a - 2} + \cdots + V^{n-2}R$ is a distributive triple
in $V^n$ for any $a$ and $n$ with $2 \leq a \leq n - 2$.

\end{corollary}

We will also need the following theorem from \cite{U}.

\begin{theorem} \label{dualkoszul}

A quadratic algebra $A$ is Koszul iff its dual algebra $A^*$ is
Koszul.  In the situation where they are both Koszul the Hilbert
series of $A$ is given by $\frac{1}{h(-x)}$ where $h(x)$ is the
Hilbert series of $A$.

\end{theorem}

\section{$Q_n$ and $Q_n(G)$}

Let $P(x) = x^n - a_{n-1}x^{n-1} + a_{n-2}x^{n-2} - \cdots +
(-1)^na_0$ be a polynomial over a division algebra.  I. Gelfand
and V. Retakh \cite{vieta} studied relationships between the
coefficients $a_i$ and a generic set $\{x_1, \cdots, x_n\}$ of
solutions of $P(x)=0$. For any ordering $(i_1, \cdots, i_n)$ of
$\{1, \cdots, n\}$ one can construct \textit{pseudoroots} $y_k$,
$k = 1, \cdots n$, (certain rational functions in $x_{i_1},
\cdots, x_{i_n}$) that give a decomposition $P(t) = (t - y_n)
\cdots (t - y_2)(t - y_1)$ where $t$ is a central variable.

In \cite{GRW} I. Gelfand, V. Retakh, and R. Wilson introduced the
algebra $Q_n$ of all pseudo-roots of a generic noncommutative
polynomial, determined a basis for this algebra and studied its
structure. The algebras $Q_n$ have a presentation given by
generators $u(A), \emptyset \neq A \subset [n]$ and relations

$$\sum_{C,D \subset A} [u(C \cup i), u(D \cup j)] = (\sum_{E \subset
A} u(E \cup i \cup j)) \sum_{F \subset A}(u(F \cup i) - u(F \cup
j))$$ for all $A \subset [n], i,j \in [n]\setminus A, i \neq j$.

In \cite{GGR} I. Gelfand, S. Gelfand, and V. Retakh introduced a
class of quotient algebras of $Q_n$ corresponding to graphs on $n$
nodes.  Let $G$ be a graph with vertex set $[n]=\{1, 2, \cdots,
n\}$ and edge set $E$ composed of elements of  $P([n])$ with
cardinality two (hence $G$ has no loops of multiple edges).  We
can then consider the quotient algebra $Q_n(G)$ we get by adding
the additional relations $u(\{i,j\}) = 0$ if $\{i,j\} \notin E$ to
$Q_n$. The following theorem gives a nice presentation of the
algebra $Q_n(G)$.

\begin{theorem} \label{graphalg} \cite{GGR} Let $G$ be a graph on $n$ nodes with edge set $E$.
Then the algebra $Q_n(G)$ is generated by the elements $u(i)$ for
$i \in [n]$ and $u(i,j)$ for $\{i,j\} \in E$ with the following
relations (assume $u(i,j)=0$ if $\{i,j\} \notin E$):

\noindent{$(i) [u(i),u(j)] = u(i,j)(u(i)-u(j)) \ i \neq j,\ i,j
\in [n]$}

\noindent{$(ii)
[u(i,k),u(j,k)]+[u(i,k),u(j)]+[u(i),u(j,k)]=u(i,j)(u(i,k)-u(j,k))$
for distinct $i,j,k \in [n]$}

\noindent{$(iii) [u(i,j),u(k,l)] = 0$ for distinct $i,j,k,l \in
[n]$}

\end{theorem}

\section{The Algebra $K_3$}

Here we consider the algebra that is generated by the graph
$K_{3}$.  By theorem \ref{graphalg} this algebra has generators
$u(1), u(2), u(3), u(12), u(13), u(23)$ together with the
following relations in $V \otimes V$ which we will refer to as
$r_{1}$ through $r_{5}$.

$r_{1} = [u(1), u(2)] + u(12)(u(2) - u(1)) = 0$

$r_{2} = [u(2), u(3)] + u(23)(u(3) - u(2)) = 0$

$r_{3} = [u(3), u(1)] + u(13)(u(1) - u(3)) = 0$

$r_{4} = [u(12), u(23)] + [u(12), u(3)] + [u(1), u(23)] -
u(13)(u(12) - u(23)) = 0$

$r_{5} = [u(12), u(13)] + [u(12), u(3)] + [u(2), u(13)] -
u(23)(u(12) - u(13)) = 0$

We have only five relations because all other possible
combinations in i) and ii) are linear combinations of these five.
We also have no relations of type iii) because $n = 3$ and we do
not have four distinct integers to work with.


We define an increasing filtration on $K_3$ by defining $F_n$ to
be the span of all monomials $u(A_1)u(A_2) \cdots u(A_k)$ such
that $ \sum_{i=1}^{k} \mid A_i \mid  \leq n$.  It is clear that
our $F_i$ are subspaces with the properties $\bigcup_i F_i = K_3$
and $F_iF_j \subseteq F_{i+j}$. We set the define $F_0$ to be the
span of 1.

Now we form $gr(K_3)$ in the usual way.  Take $G_i$ = $F_i /
F_{i-1}$ and set $gr(K_3)$ = $\bigoplus_i G_i$ and then define
multiplication in $gr(K_3)$ so for all $a \in F_i, b \in F_j, (a +
F_{i-1})(b + F_{j-1}) = ab + F_{i+j-1}$.  Note that there is a
non-linear map $gr: K_3 \rightarrow gr K_3$ that sends $a \in F_i,
 a \notin F_{i-1}$ to $a + F_{i-1}$ in $gr(K_3)$ and sends 0 to 0.

\section{A New Presentation of $gr(K_3)$}

For ease of notation, let us temporarily set $a=u(1), b=u(2),
c=u(3), d=u(12), e=u(23),$ and $f=u(13)$.  Notice our relations in
$K_3$ then become:

$r_1 = db - da + ab - ba$

$r_2 = ec - eb + bc - cb$

$r_3 = fa - fc + ca - ac$

$r_4 = de - ed - fd + fe + dc - cd + ae - ea$

$r_5 = df - fd - ed + ef + dc - cd + bf - fb$

Therefore in $gr(K_3)$ (if we allow each generator to represent
itself under the image $gr$) we know the following relations hold:

$d \cdot b = da$

$e \cdot c = eb$

$f \cdot c = fa$

$f \cdot e = fd - de + ed$

$f \cdot d = df + ef - ed$

Now this list of relations might not be enough for a presentation
of $gr(K_3)$, for other relations may hold true and be needed as
well.  Call the algebra generated by the five truncated relations
above ``chopped'' $K_3$ (or $ch(K_3)$).  We know that since
$gr(K_3)$ is a quotient of $ch(K_3)$ (it has possibly more
relations) we can verify that the two are equal by showing they
have the same Hilbert series.  Since $gr(K_3)$ has the same
Hilbert series as $K_3$ we can compare the series for $ch(K_3)$ to
the one for $K_3$ instead.

First let us use the diamond lemma to compute the Hilbert series
of $K_3$.

\begin{theorem} \label{K3basis}

A basis for $K_3$ is given by the set of monomials in $T(V)$
containing none of the following substrings: $cef, cd, cb, ca, bf,
ba$.

\end{theorem}

\begin{proof}

Order the set of monomials first by monomial length, and then
lexicographically with the ordering $c > b > e > f > a > d$ (which
is actually $u(3) > u(2) > u(23) > u(13) > u(1) > u(12)$).  We
then have the following reductions (after replacing $r_5$ with
$r_5 - r_4$):

$cd \rightarrow de - ed - fd + fe + dc + ae - ea$

$bf \rightarrow fb + ae - ea + fe - ef + de - df$

$ba \rightarrow ab + db - da$

$cb \rightarrow bc + ec - eb$

$ca \rightarrow ac + fc - fa$

This gives us $cba$ and $cbf$ as two ambiguities that need to be
resolved.  It is not hard to check that if we compute $c(ba) -
(cb)a$ and reduce with the above five relations we get zero.
However in order to resolve $cbf$ it turns out we must add the
relation: $cef = cfb + cfe + ace + fce - fae - cea + dce + de^2 +
ae^2 + fe^2 - fde - dcf - def - aef + eaf - fef + fdf - bcf - ecf
+ efb + efe - e^2f - e^2a - edf$

However, adding this relation creates no new ambiguities.  So we
have found a basis in the set of monomials not containing the
strings $cef, cd, cb, ca, bf, ba$.

\end{proof}

\begin{theorem}

The Hilbert series of $K_3$ is $H(x) = \frac{\rm 1}{\rm x^3 - 6x^2
+ 5x - 1}$.

\end{theorem}

\begin{proof}

We must count the number of monomials of length $n$ in $T(V)$ not
containing any of the strings listed in theorem \ref{K3basis}.
Call such monomials the \textit{valid} monomials of length $n$.

Let $T_n$ be the number of valid monomials of length $n$.

Let $J_n$ be the number of valid monomials of length $n$ that
begin with $b$.

Let $K_n$ be the number of valid monomials of length $n$ that
begin with $c$.

We can note right away that $T_{n+1} = 4 T_n + J_{n+1} + K_{n+1}$.
Since words beginning with $b$ can not be followed by $a$ or $f$,
we get $J_{n+1} = 2 T_{n-1} + J_n + K_n$.  Words beginning with
$c$ can be followed by $c, f,$ or $e$, but in the e case they can
not next be followed by $f$.  This gives us $K_{n+1} = T_{n-1} +
K_n + 3 T_{n-2} + J_{n-1} + K_{n-1}$.  By counting the valid words
up to length three, we also get initial conditions.  Thus we
obtain the following system of recurrences:

$T_{n+1} = 4 T_n + J_{n+1} + K_{n+1}$

$J_{n+1} = 2 T_{n-1} + J_n + K_n$

$K_{n+1} = T_{n-1} + K_n + 3 T_{n-2} + J_{n-1} + K_{n-1}$

$T_1 = 6, T_2 = 31, T_3 = 157$

To solve this system notice first that in the second equation $J_n
+ K_n = J_{n+1} - 2 T_{n-1}$.  Plugging this into our first
equation gives $J_{n+2} = T_{n+1} - 2 T_n.$  We can use this to
get rid of the $J$'s in the first and third equations to get the
system:

$T_{n+1} = K_{n+1} + 5 T_n - 2 T_{n-1}$

$K_{n+1} = K_n + K_{n-1} + T_{n-1} + 4 T_{n-2} - 2 T_{n-3}$

Solving the first for $K_{n+1}$ and substituting into the second
gives us $T_{n+1} = 6T_n - 5T_{n-1} + T_{n-2}$.  Using generating
functions and our initial conditions we can quickly find that the
Hilbert series for $K_3$ is $H(x) = \frac{\rm 1}{\rm x^3 - 6x^2 +
5x - 1}$.

\end{proof}

Now we must find the Hilbert series of $ch(K_3)$.  Consider an
ordering first by monomial length and then lexicographically with
$f > e > d > c > b > a$.  We get the following reductions in
$ch(K_3)$.

$fe \rightarrow ef + df - de$

$fd \rightarrow df + ef - ed$

$db \rightarrow da$

$ec \rightarrow eb$

$fc \rightarrow fa$

We have two ambiguities to resolve this time: $fec$ and $fdb$.  An
attempt to resolve the $fec$ ambiguity shows it necessary to add
the relation $efb = efa - dfb + dfa$.  With this new relation, we
can resolve the $fdb$ ambiguity.  However, we created a new
ambiguity by adding our $efb$ relation.  In order to resolve
$fefb$ we must toss in the relation $effb = effa + dffa - dffb +
\frac{1}{2}edfb - \frac{1}{2}edfa - \frac{1}{2}ddfb + \frac{1}{2}
dd fa$.  We now have to worry about the ambiguity $feffb$.  We can
deal with all these ambiguities at once with the following lemma.

\begin{lemma}

Suppose we need to resolve an ambiguity of the form $ef^nb = ev_n
+ dw_n - \alpha_n df^nb$ where $v_n$ and $w_n$ are linear
combinations of monomials of length $n+1$. Suppose also that the
terms in $v_n$ and $w_n$ are all less than or equal then $f^{n}b$
and that $\alpha_n$ is a positive real number.  Then the ambiguity
$fef^nb$ can be resolved by adding a relation of the form
$ef^{n+1}b = ev_{n+1} + dw_{n+1} - \alpha_{n+1}df^{n+1}b$ where
$v_{n+1}$ and $w_{n+1}$ are linear combinations of monomials of
length $n+2$, the terms in $v_{n+1}$ and $w_{n+1}$ are all less
than $f^{n+1}b$, and $\alpha_{n+1}$ is a positive real number.

\end{lemma}

\begin{proof}

We have $(fe)f^nb - f(ef^nb) = ef^{n+1}b + df^{n+1}b - (d+f)(ev_n
+ dw_n - \alpha_ndf^nb) = ef^{n+1}b + df^{n+1}b - dev_n - d^2w_n -
\alpha_nd^2f^nb - efv_n - dfv_n + dev_n - dfw_n -efw_n +
\alpha_ndf^{n+1}b + \alpha_nef^{n+1}b - \alpha_nedf^{n}b$ .

Thus $ef^{n+1}b + \alpha_nf^{n+1}b = d(-f^{n+1}b + dw_n +
\alpha_ndf^nb + fv_n + fw_n - \alpha_nf^{n+1}b) + e(fv_n + fw_n -
dw_n + \alpha_ndf^nb)$ and dividing by $1 + \alpha_n$ gives us:

$ef^{n+1}b = \frac{1}{1 + \alpha_n}e(fv_n + fw_n - dw_n +
\alpha_ndf^nb) + \frac{1}{1 + \alpha_n}d(-f^{n+1}b + dw_n +
\alpha_ndf^nb + fv_n + fw_n) - \frac{\alpha_n}{1 +
\alpha_n}df^{n+1}b$

Setting $v_n = \frac{1}{1 + \alpha_n}(fv_n + fw_n - dw_n +
\alpha_ndf^nb), w_n = \frac{1}{1 + \alpha_n}(-f^{n+1}b + dw_n +
\alpha_ndf^nb + fv_n + fw_n)$ and taking $\alpha_{n+1}$ to the
positive real constant $\frac{\alpha_n}{1 + \alpha_n}$ completes
the proof.

\end{proof}

With this lemma, we see that the only bad words are ones
containing strings of the following forms: $fe, fd, db, ec, fc,$
or $ef^nb$ for $n \geq 1$.  We can now use this to find the
Hilbert series of $ch(K_3)$.

\begin{proposition}

The Hilbert series of $ch(K_3)$ is equal to the Hilbert series of
$K_3$.

\end{proposition}

\begin{proof}

We wish to count the strings not containing $fe, fd, db, ec, fc,$
or $ef^nb$ for $n \geq 1$ as a substring.  To count all such
strings let $T_n$ be the total number of valid monomials of length
$n$. Let $K_n$ be the number of valid monomials of length $n$
beginning with $d$. Let $L_n$ and $M_n$ be the corresponding
numbers for $f$ and $e$ respectively.

We can immediately see that $T_n = K_n + L_n + M_n + 3 T_{n-1}$.
We know that if a word begins with $d$, it can be followed by any
smaller valid word not beginning with $b$.  This gives us $K_n =
2T_{n-2} + K_{n-1} + L_{n-1} + M_{n-1}$.  The words beginning with
$f$ can be followed by $f, a$ or $b$ giving us $L_n = 2T_{n-2} +
L_{n-1}$.

Finally we must count the words beginning with $e$.  If $e$ is
followed by any valid word not beginning with $f$ or $c$ then we
are okay.  There will be $2 T_{n-2} + K_{n-1} + M_{n-1}$ of these.
If the second letter does happen to be $f$ then the next letter
can only be $a$ or $f$.  If it is $a$ then we can follow up with
any valid word (which adds $T_{n-3}$ to the equation), but if it
is $f$ we are once again in an $a$ or $f$ situation.  This time
the $a$ case ends up adding $T_{n-4}$.  We can repeat this down
the line to add $T_{n-5} + T_{n-6} + \cdots + T_1$ and finally we
add 2 (or $2T_0$) for the $ef \cdots fa,$ and $ef \cdots ff$
cases. This gives us $M_n = 2T_{n-2} + K_{n-1} + M_{n-1} +
\sum_{k=3}^n T_{n-k} + T_0$.

We must now solve the following system of recurrences:

$T_n = K_n + L_n + M_n + 3 T_{n-1}$

$K_n = 2T_{n-2} + K_{n-1} + L_{n-1} + M_{n-1}$

$L_n = 2T_{n-2} + L_{n-1}$

$M_n = 2T_{n-2} + K_{n-1} + M_{n-1} + \sum_{k=3}^n T_{n-k} + T_0$

To do this, first we must get rid of the summation for $M_n$.  Set
$R_n = M_n - M_{n-1}$ which equals $2 T_{n-2} + K_{n-1} + M_{n-1}
+ \sum_{k=3}^n T_{n-k} + T_0 - 2 T_{n-3} - K_{n-2} - M_{n-2} -
\sum_{k=3}^{n-1} T_{n - 1 - k} - T_0 = 2 T_{n-2} - 2 T_{n-3} +
K_{n - 1} - K_{n-2} + M_{n-1} - M_{n-2} + \sum_{k=3}^n T_{n-k} -
\sum_{k=4}^{n} T_{n - k} = 2 T_{n-2} - 2 T_{n-3} + K_{n - 1} -
K_{n-2} + R_{n-1} + T_{n-3}$ so $R_n = 2 T_{n-2} - T_{n-3} + K_{n
- 1} - K_{n-2} + R_{n-1}$.

Now we can get rid of our $M_n$ terms altogether by replacing
$T_n$ and $K_n$ with $T_n - K_n = 3T_{n-1} - 2 T_{n-2} + K_n -
K_{n-1} + L_n - L_{n-1} + R_n$ and $T_n - K_{n+1} = 3T_{n-1} + K_n
+ L_n + M_n - 2T_{n-1} - K_n - L_n - M_n = T_{n-1}$.

By plugging $L_n - L_{n-1}= 2T_{n-2}$ into our $T_n$ equation we
are left with the following system:

$T_n = 3T_{n-1} + 2 K_n - K_{n-1} + R_n$

$R_n = 2T_{n-2} - T_{n-3} + R_{n-1} + K_{n-1} - K_{n-2}$

$K_n = T_{n-1} - T_{n-2}$

Substituting for $K_n$ leads to the system:

$T_n = 5T_{n-1} - 3T_{n-2} + T_{n-3} + R_n$

$R_n = 3T_{n-2} - 3T_{n-3} + T_{n-4} + R_{n-1}$

We can now solve for $R_n$ in the first equation and substitute
into the second equation to get $T_n = 6T_{n-1} - 5T_{n-2} +
T_{n-3}$.  This is the same recurrence we used to generate the
Hilbert series of $K_3$.  Checking the length one, two, and three
cases gives the same initial conditions as well.  Therefore the
two algebras must have the same Hilbert series.

\end{proof}

\begin{corollary}

$ch(K_3) \cong gr(K_3)$.

\end{corollary}

\begin{proof}

We know $ch(K_3)$ has the same graded dimension as $K_3$ which has
the same graded dimension as $gr(K_3)$.  Since $ch(K_3)$ is a
quotient of $gr(K_3)$ with the same graded dimension, the two must
be isomorphic.

\end{proof}

From now on we will list the chopped relations for a presentation
of $gr(K_3)$.

\section{Reduction to $gr(K_3)$}

The following corollary from \cite{PP} will be useful here.

\begin{corollary}

Let $0 = F_0W \subset F_1W \subset \cdots \subset F_{l-1}W \subset
F_lW$ be a filtered vector space and $X_1, \cdots, X_n$ be a
collection of subspaces; then the following two conditions are
equivalent:

(a) the whole set of subspaces $F_0W, \cdots, F_lW, X_1, \cdots,
X_n \subset W$ is distributive.

(b) the associated graded collection $gr X_1, \cdots, gr X_n$ in
the associated graded vector space $gr W$ is distributive are for
any $1 \leq i < j \leq n$ either of the two equivalent conditions
holds:

$gr (X_i + X_j) = gr X_1 + gr X_j$ or $gr(X_i \cap X_j) = gr X_i
\cap gr X_j$.

\end{corollary}

Hence, if the set $\{gr X_i\}_i$ generates a distributive lattice
in $gr V$, and $gr(X_i \cap X_j) = gr(X_i) \cap gr (X_j)$ for all
$i$ and $j$, then the set $\{X_i\}_i$ generates a distributive
lattice in $V$. We wish to check that $gr(RV \cap VR) = gr(RV)
\cap gr(VR)$. By looking at $gr$ as a function and only using the
basic rules for set maps and intersections we see that $gr(RV \cap
VR) \subseteq gr(RV) \cap gr(VR)$.  It will be enough for us to
show that the dimension of $gr(RV) \cap gr(RV)$ is one and that
$gr(RV \cap VR)$ is not the zero subspace.

For the second part, notice that the vector $r_1c + r_2a + r_3b +
r_4(c-b) + r_5(a-c)$ is equal to the vector $ar_2 + br_3 + cr_1 +
d(r_2 + r_3) + e(r_1 + r_3) + f(r_1 + r_2)$ where $r_1, \cdots,
r_5$ are the relations we defined earlier for $K_3$.  This shows
that $RV \cap VR$ is not zero, and thus $gr(RV \cap VR)$ is not
the zero subspace.

\begin{proposition}

dim $gr(RV) \cap gr(VR) = 1$

\end{proposition}

\begin{proof}

In the last section we showed that $gr(R)$ is the span of $\{fe -
ef - df + de, fd - df - ef + ed, d(b-a), e(c-b), f(a-c) \}$. Since
these are the relations we will be working with throughout the
rest of this paper, we officially set:

$r_1 = db - da$

$r_2 = ec - eb$

$r_3 = fa - fc$

$r_4 = de - ed - fd + fe$

$r_5 = df - fd - ed + ef$

Call the span of these five relations $S$ (so $S = gr(R)$).  Our
goal is to show $SV \cap VS$ is of dimension one.

Suppose $x \in SV \cap VS$ and $x \not= 0$.  Since all the
monomials in $VS$ contain no $a, b,$ or $c$ in the middle spot, we
can replace $SV$ with $sp\{r_4, r_5\} \otimes V$.  As all the
monomials in $SV$ contain no $a, b$, or $c$ in the first slot we
can replace $VS$ with $sp\{d,e,f\} \otimes S$.

Suppose $x \in sp\{r_4, r_5\} \otimes V \cap sp\{d,e,f\} \otimes
S$.  Then we can write $x = r_4v_1 + r_5v_2$ for some $v_1,v_2 \in
V$.  Now if the coefficient of $c$ in $v_1$ was nonzero we would
have an $fee$ term with nothing else that could cancel it out.
Hence we would have an $fee$ appearing in $VS$ which is not
possible. If the coefficient of $f$ in $v_1$ was nonzero then we
would have a $def$ which could only happen in $VS$ is the
coefficient of $dr_5$ was nonzero.  But this would give us a
nonzero $dee$ term with nothing to cancel it out, which is also
not possible.  If the coefficient of $d$ was nonzero we would have
an $fed$ appearing which implies that $\alpha fr_4 + \beta fr_5$
appears in VS with $\alpha \not= - \beta$.  As $ffe$ can not
appear, $\alpha$ must be zero, so $\beta \not= 0$.  But then we
would have an $-\beta$ $ffd$ appearing with nothing to cancel it
out with.  This shows that $v_1 \in sp\{a,b,c\}$.  A similar
argument shows that $v_2$ is in this same span.  We know now that
$x \in SV \cap VS$ implies $x \in sp\{r_4, r_5\} \otimes
sp\{a,b,c\} \cap sp\{d,e,f\} \otimes S$.  We can replace this last
$S$ with $sp\{r_1, r_2, r_3\}$ after noticing that only $a$, $b$,
and $c$ can now appear in the last slot.  So far we know $SV \cap
VS$ = $sp\{r_4, r_5\} \otimes sp\{a,b,c\} \cap sp\{d,e,f\} \otimes
sp\{r_1, r_2, r_3\}$

Now suppose $x = r_4v_1 + r_5v_2$ where $v_1,v_2 \in sp\{a,b,c\}$
and $x \not= 0$.  Suppose also that the coefficient of $b$ in
$v_1$ is 0.  Then no $deb$ or $feb$ can appear in $x$ (with a
non-zero coefficient).  Looking in $VS$ we see this means $dr_2$
and $fr_2$ must have coefficients of zero.  This means no $dec$ or
$fec$ can appear in $x$ either.  Hence $v_1$ is a constant
multiple of $a$. However this constant must be 0, otherwise the
term $dea$ would appear, and there is no way to achieve that in
$VS$.  So in this case $v_1$ is 0 and $x = r_5v_2$.  Then $v_2$
would have to be a multiple of $a$ since no $dfb$ and $fdc$ can
occur.  So $x$ is a constant multiple of $r_5 a$.  But $r_5 a$ is
not in $VS$ so we have reached a contradiction.

We now know that if $x = r_4v_1 + r_5v_2 \in SV \cap VS$ then the
coefficient of $b$ in $v_1$ is non-zero, so we can scale $x$ to
make this coefficient 1.  Hence $deb + feb$ appears in $x$.
Looking in $SV$ we see this implies $-dec -fec$ appears as well.
This implies the coefficient of $c$ in $v_1$ is $-1$.  As $fdc$
can not appear, this term must cancel, meaning the coefficient of
$c$ in $v_2$ is 1. This means $dfc + efc$ appears in $x$.  Looking
in $SV$ we see $-dfa -efa$ appears, implying the coefficient of
$a$ in $v_2$ is $-1$.  Since the coefficients of $a$ in $v_1$ and
$b$ in $v_2$ must be 0 (look at $SV$ to see this) we get that $x =
r_4(b-c) + r_5(c-a)$. It is easy to see this $x$ is in $VS$
because it is equal to $(e+f)r_1 - (d+f)r_2 - (d+e)r_3$.

\end{proof}

\section{A Reduction to Two Cases}

Recall from corollary \ref{mainkoszul} that $K_3$ will be Koszul
if we show that $RV^{n - 2} \cap VRV^{n - 2} \cap \cdots \cap V^{a
- 2}RV^{n - a}$, $V^{a - 1}RV^{n - a - 1}$, $V^{a}RV^{n - a - 2} +
\cdots + V^{n-2}R$ is a distributive triple in $V^n$ for any $n$
and $2 \leq a \leq n - 2$.

We will require the following lemma of Serconek and Wilson (lemma
1.1 from \cite{qnkoszul}):

\begin{lemma} \label{SWchop}

If

1)  $V = \sum_{i \in I} V_i$ is graded as a vector space

2)  $X_j$ is a collection of subspaces of $V$

3)  Each $X_j = \sum_{i \in I} (X_j \cap V_i)$

\noindent then $\{X_j\}_j$ is distributive if and only if for all
$i \in I$, $\{X_j \cap V_i\}_j$ is distributive in $V_i$.

\end{lemma}

We now choose a particular $\{1,2\}^n$ grading of $V$ to apply
this lemma to.  Set $V^n_{(i_1, i_2, \cdots, i_n)}$ to be the span
of all monomials $u(A_1)u(A_2) \cdots u(A_n)$ so that $|A_k| =
i_k$. For example, in the n = 2 case $r_1$, $r_2$, and $r_3$ are
in the $V^2_{(2,1)}$ space and $r_4$ and $r_5$ are in the
$V^2_{(2,2)}$ space.  Hence $R = (R \cap V^2_{(2,1)}) + (R \cap
V^2_{(2,2)})$. This shows that property three of lemma
\ref{SWchop} will apply to our sets $\{V^{a-2}RV^{n-a}\}$.

We have to show that $\{RV^{n-2} \cap (V^n)_{\alpha}, \cdots,
V^{n-2}R \cap (V^n)_{\alpha} \}$ generates a distributive lattice
in $(V^n)_\alpha$ for $\alpha \in \{1,2\}^n$.  Notice that we need
only check the $\alpha = (i_1, i_2, \cdots, i_n)$ such that $i_1
\geq i_2 \geq \cdots \geq i_n$.  This is because $R \in
V^2_{(2,1)} + V^2_{(2,2)}$ so if a $(1,2)$ appears somewhere in
the string $\alpha$ then one of our $V^{a-2}RV^{n-a} \cap
(V^n)_{\alpha}$ will be zero.  Since we can show proper subsets of
$\{RV^{n-2} \cap (V^n)_{\alpha}, \cdots, V^{n-2}R \cap
(V^n)_{\alpha} \}$ are distributive, we are done for such
$\alpha$.

Next notice that the subspace $RV \cap VR$ we found earlier is
contained in $V^3_{(2,2,1)}$.  Combining corollary
\ref{mainkoszul} with lemma \ref{SWchop}, we have to check that
$RV^{n-2} \cap \cdots \cap V^{a-2}RV^{n-a} \cap (V^n)_{\alpha},
V^{a-1}RV^{n-a-1} \cap (V^n)_{\alpha}, (V^{\alpha}RV^{n-a-2} \cap
(V^n)_{\alpha})+ \cdots + (V^{n-2}R \cap (V^n)_{\alpha})$ is a
distributive triple for any decreasing $\alpha$ and $2 \leq a \leq
n-2$.  However if the last two digits of $\alpha$ are (1,1) then
the last term of the third element $(V^{n-2}R \cap (V^n)_{\alpha})
= 0$ and we will be done because proper subsets are distributive.
We are done with all cases except when $\alpha$ contains a two in
the second to last spot. This leaves $\alpha = (2,2,\cdots,2,1)$
and $\alpha = (2,2,\cdots,2)$.

Next suppose $a > 2$.  Notice that our first term in our triple is
$RV^{n-2} \cap VRV^{n-3} \cap \cdots \cap V^{a-2}RV^{n-a} \cap
(V^n)_{\alpha}$ which is contained in $(RV \cap VR)V^{n-3} \cap
(V^n)_{\alpha}$.  Since $\alpha$ can not start out with $(2,2,1)$
and $RV \cap VR$ is contained in this graded space, this term must
be zero.  Hence we need only check the case where $a = 2$.

From here on when we say a set $\{X_1, \cdots, X_n\}$ is
distributive in the $\alpha$ case we mean that the set $\{X_1 \cap
(V^n)_{\alpha}, \cdots, X_n (V^n)_{\alpha}\}$ is distributive.
Thus we must show $RV^{n-2},$ $VRV^{n-3},$ $V^2RV^{n-4} + \cdots +
V^{n-2}R$ is a distributive triple in the $(2, 2, \cdots, 2)$ and
$(2, 2, \cdots, 2, 1)$ cases.  This amounts to showing $RV^{n-2}
\cap (VRV^{n-3} + V^2RV^{n-4} + \cdots + V^{n-2}R) = (RV^{n-2}
\cap VRV^{n-3}) + (RV^{n-2} \cap (V^2RV^{n-4} + \cdots +
V^{n-2}R))$ for those two cases.  But as $RV \cap VR$ is in
$V^3_{(2,2,1)}$, we know that $RV^{n-2} \cap VRV^{n-3} = 0$ in
these cases.  This simplifies what we must show to $RV^{n-2} \cap
(VRV^{n-3} + V^2RV^{n-4} + \cdots + V^{n-2}R) = (RV^{n-2} \cap
(V^2RV^{n-4} + \cdots + V^{n-2}R))$. And since the right side is
contained in the left we only have one direction left to show.  We
need show that if $x \in RV^{n-2}$ and $x \in VRV^{n-3} + \cdots +
V^{n-2}$ then $x \in V^2RV^{n-4} + \cdots + V^{n-2}R$.

Introducing a little more terminology makes this statement
simpler.  Define $W_k = RV^{k-2} + VRV^{k-3} + \cdots + V^{k-2}R$
for $k \geq 2$ and set $W_k$ to the zero subspace otherwise.  Then
we must show that $X \in RV^{n-2} \cap VW_{n-1} \Longrightarrow x
\in V^2W_{n-2}$ for our two cases.  Notice also that this
statement is trivial for $n < 4$.

Next we look for a more natural spanning set for $R$.  Let
$\sigma$ be the permutation (1,2,3).  Notice that the map $T$
sending $u(A)$ to $u(\sigma(A))$ sends $sp\{r_1, r_2, r_3\}$ and
$sp\{r_4, r_5\}$ both back to themselves.  Hence $R$ is $T$
invariant. Letting $\omega$ be a primitive cube root of one, the
following elements form a basis for $R$ consisting only of
eigenvectors. Set:

$u_1 = u(12) + u(23) + u(13)$

$u_\omega = u(12) + \omega u(23) + \omega^2u(13)$

$u_{\omega^2} = u(12) + \omega^2u(23) + \omega u(13)$

$v_1 = u(1) + u(2) + u(3)$

$v_2 = u(1) + \omega u(2) + \omega^2 u(3)$

$v_3 = u(1) + \omega^2 u(2) + \omega u(3)$

\noindent then our relations become

$r_1 = (u_1 + u_\omega + u_{\omega^2})(\omega - 1)v_1 + (\omega^2
- 1)v_2)$

$r_2 = (\omega^2u_1 + u_\omega + \omega u_{\omega^2})((\omega^2 -
\omega)v_1 + (\omega - \omega^2)v_2)$

$r_3 = (\omega u_1 + u_\omega + \omega^2u_{\omega^2})((1 -
\omega^2)v_1 + (1 - \omega)v_2)$

$r_4 = u_{\omega^2}^2 - 2u_1u_\omega + u_\omega u_1$

$r_5 = u_\omega^2 - 2u_1u_{\omega^2} + u_{\omega^2}u_1$

If we set $a = v_1, b = v_2, d = u_1, e = u_\omega, f =
u_{\omega^2}$ then we get the simpler looking set:

$r_1 = (e + f + d)((\omega - 1)a + (\omega^2 - 1)b)$

$r_2 = (e + \omega f + \omega^2 d)((\omega^2 - \omega)a + (\omega
- \omega^2)b)$

$r_3 = (e + \omega^2 f + \omega d)((1 - \omega^2)a  + (1 -
\omega)b)$

$r_4 = f^2 - 2de - ed$

$r_5 = e^2 - 2df + fd$

Keep in mind that $r_4$ and $r_5$ still sit in $V^2_{(2,2)}$ and
$r_1, r_2,$ and $r_3$ still sit in $V^2_{(2,1)}$ since $d$, $e$,
and $f$ are multiples of $u(12)$, $u(23)$, and $u(12)$ and $a$,
$b$, and $c$ are multiples of $u(1)$, $u(2)$, and $u(3)$.  This
will be our spanning set for $R$ as we move on to the two last
cases.

\section{The Two Last Cases}

Our main goal is to prove the following lemma in the
$(2,2,\cdots,2)$ and $(2,2,\cdots,1)$ cases.

\begin{lemma} \label{lastcaselemma}

Suppose $n \geq 2$ then

a)  If $z_1, z_2 \in V^n, dz_1 + ez_2 \in W_{n+1}$ then $z_1, z_2
\in W_n$

b)  If $y_1, y_2 \in V^n, ey_1 + fy_2 \in W_{n+1}$ then $y_1, y_2
\in W_n$

c)  If $x_1, x_2 \in V^n, dx_1 + fx_2 \in W_{n+1}$ then $x_1, x_2
\in W_n$

\end{lemma}

\begin{proof}

Suppose we knew the lemma was true for either $n = 1$ or $n = 2$.
Then we can assume by induction that the lemma holds for $n - 1$.
Suppose we are in the a) case and set $z = dz_1 + ez_2 \in
W_{n+1}$.  Then we can write $z = r_4h_1 + r_5h_2 + VW_n = f(fh_1
+ dh_2) + e(eh_2 + dh_1) - 2d(eh_1 + fh_2) + VW_n$.  This means
$fh_1 + dh_2 \in W_n$.  Since $z \in dV^n + eV^n$ we know that
$f(fh_1 + dh_2) \in VW_n$ so $fh_1 + dh_2 \in W_n$.  By our
inductive hypothesis $h_1, h_2 \in W_{n-1}$.  So $z \in r_4W_{n-1}
+ r_5W_{n-1} + VW_n \subset VW_n$ and hence $z_1, z_2 \in W_n$.
The b) and c) cases are similar.

We have left to find basis cases for lemma.  In the
$(2,2,\cdots,2)$ situation we can find one when $n=1$.  Assume
$dz_1 + ez_2 \in W_2 = R = \alpha r_4 + \beta r_5$.  If $\alpha
\not= 0$ we would have an $f^2$ appearing which we could not
cancel, and hence a contradiction.  Similarly $\beta$ must be 0 as
well and hence $z_1$ and $z_2$ are both in $W_1 = 0$.

\end{proof}

Now that we know that the lemma is true in both cases we can prove
the following proposition thus completing our proof that $K_3$ is
Koszul.

\begin{proposition}

If $n \geq 4, x \in RV^{n-2} \cap V_{n-1}$ then $x \in V^2W_{n-2}$

\end{proposition}

\begin{proof}

Since $x \in RV^{n-2}$ we can write $x = r_4h_1 + r_5h_2 = (f^2 -
2de - ed)h_1 + (e^2 - 2df + fd)h_2 = f(fh_1 + dh_2) + e(eh_2 +
dh_1) - 2d(eh_1 + fh_2)$.  Since $x \in VW_{n-1}$ we know that
$fh_1 + dh_2, eh_2 + dh_1$, and $eh_1 + fh_2$ are all in
$W_{n-1}$.  Since $n-1 \geq 3$, lemma \ref{lastcaselemma} applies
and thus $h_1, h_2 \in W_{n-2}$.  Thus $x \in r_4W_{n-2} +
r_5W_{n-2} \subset V^2W_{n-2}$ and we are done.

\end{proof}

\begin{theorem}

$K_3$ is Koszul.

\end{theorem}

\bibliographystyle{plain}

\bibliography{citations}

\end{document}